\documentclass[12pt]{amsart}

\usepackage[section]{placeins}

\usepackage[utf8]{inputenc}  

\usepackage{fullpage}

\usepackage{amsmath}
\usepackage{amsfonts}
\usepackage{amssymb}
\usepackage{amsthm}
\usepackage{latexsym}

\usepackage{enumerate}
\usepackage{cancel}
\usepackage{cases}
\usepackage{empheq}
\usepackage{multicol}

\usepackage{graphicx} 

\theoremstyle{plain}
\newtheorem{theorem}{Theorem}[section]

\theoremstyle{definition}

\theoremstyle{remark}
\newtheorem{remark}[theorem]{Remark}

\numberwithin{equation}{section} 
\numberwithin{figure}{section}   

\newcommand{\norm}[1]{\left\|#1\right\|}

\title[Data Assimilation with Moving Assimilation Points]{Continuous Data Assimilation with a Moving Cluster of Data Points for a Reaction Diffusion Equation: A Computational Study}

\date{\today}

\author{Adam Larios}
\address[Adam Larios]{Department of Mathematics, 
	University of Nebraska--Lincoln,
	Lincoln, NE 68588-0130, USA}
\email[Adam Larios]{alarios@unl.edu}
\author{Collin Victor}
\address[Collin Victor]{Department of Mathematics, 
	University of Nebraska--Lincoln,
	Lincoln, NE 68588-0130, USA}
\email[Collin Victor]{collin.victor@huskers.unl.edu}
\thanks{MSC 2010 Classification: 
	35K57 
	35K40, 
	35K61, 
	37C50, 
	35Q93, 
	34D06.}

\usepackage{caption}

\begin{document}

\begin{abstract}
Data assimilation is a technique for increasing the accuracy of simulations of solutions to partial differential equations by incorporating observable data into the solution as time evolves.  Recently, a promising new algorithm for data assimilation based on feedback-control at the PDE level has been proposed in the pioneering work of Azouani, Olson, and Titi (2014).  The standard version of this algorithm is based on measurement from data points that are fixed in space.  In this work, we consider the scenario in which the data collection points move in space over time.  We demonstrate computationally that, at least in the setting of the 1D Allen-Cahn reaction diffusion equations, the algorithm converges with significantly fewer measurement points, up to an order or magnitude in some cases.  We also provide an application of the algorithm to an inverse problem in the case of a uniform static grid.

\end{abstract}


\keywords{Continuous Data Assimilation, Allen-Cahn, Reaction-Diffusion, Moving Mesh}

\maketitle

\date{\today}

\section{Introduction}\label{secInt}
\noindent
Recently, a new method has emerged as a promising approach to data assimilation. 
We refer to this method as the Azouani-Olson-Titi (AOT) algorithm after the authors who pioneered this idea in \cite{Azouani_Olson_Titi_2014,Azouani_Titi_2014}.  
In the standard implementation of the AOT algorithm, data is collected at points on a static grid.
Here, we examine the effect of allowing the data points to move in time.
As we demonstrate below, this can lead to order-of-magnitude improvements on rates of convergence and on the number of data points needed for convergence to occur.  
Indeed, order of magnitude improvements were observed over a wide range of parameters, including those in which the equation is more computationally demanding.
In this study, we focus on the particular case of a localized cluster of data points moving at a uniform speed throughout the domain.  
Such a scenario may arise in realistic settings; e.g., a moving probe in an experiment, a vehicle or drone loaded with sensors as it crosses a crop, or a satellite sampling data as it orbits.   
Rather than study such complex settings however, we examine this technique in the context of a relatively simple equation, namely the one-dimensional Allen-Cahn reaction-diffusion equation.  

In addition, we consider an application to an inverse problem. Namely, we use the AOT algorithm to find a relationship between the viscosity and the minimum wave length of a transition layer. We experimentally estimate the number of data assimilation nodes for the uniform static grid required for convergence using periodic initial data. 
In particular, we find a relationship between the minimum required number of nodes and the minimum length scale for the Allen-Cahn equation using AOT data assimilation combined with a statistical technique we develop below. 

The term \textit{data assimilation} refers to a class of schemes that employ observational data in simulations.  
It is the subject of a large body of work (see, e.g., \cite{Daley_1993_atmospheric_book,Kalnay_2003_DA_book,Law_Stuart_Zygalakis_2015_book}, and the references therein).   
By using incoming data, data assimilation techniques aim to increase the accuracy of solutions and obtain better estimates of initial conditions. 
This incoming data is used in simulations to drive the system to the ``true'' solution. 
This strategy is widely used in many simulation-driven fields where it is almost impossible to have complete initial data, such as numerical prediction of weather on earth and on the sun, the evolution of soil moisture, biophysical monitoring, and many other areas. 

Classical data assimilation techniques are based on the Kalman Filter, a form of linear quadratic estimation.  
There are also variational methods collectively known as 3D/4D Var techniques.  
These methods are described in detail in several textbooks, including \cite{Daley_1993_atmospheric_book,Kalnay_2003_DA_book,Law_Stuart_Zygalakis_2015_book}, and the references therein. 
The AOT algorithm is an entirely different approach that adds a feedback control term at the PDE level. 
A similar approach is followed in \cite{Blomker_Law_Stuart_Zygalakis_2013_NL} in the context of stochastic differential equations.  We note that the AOT method superficially appears similar to the nudging or Newtonian relaxation methods introduced in \cite{Anthes_1974_JAS,Hoke_Anthes_1976_MWR}; however, the use of the interpolation in the AOT method is a major difference between the two methods, with crucial effects in terms of implementation, convergence rates, and the amount of measurement data required.  For an overview of nudging methods, see, e.g., \cite{Lakshmivarahan_Lewis_2013} and the references therein.

A large amount of recent literature has built upon the AOT algorithm; see, e.g., \cite{Albanez_Nussenzveig_Lopes_Titi_2016,Altaf_Titi_Knio_Zhao_Mc_Cabe_Hoteit_2015,Bessaih_Olson_Titi_2015,Biswas_Foias_Monaini_Titi_2018downscaling,Biswas_Hudson_Larios_Pei_2017,Biswas_Martinez_2017,Carlson_Hudson_Larios_2018,Celik_Olson_Titi_2018,Farhat_GlattHoltz_Martinez_McQuarrie_Whitehead_2018,Farhat_Jolly_Titi_2015,Farhat_Lunasin_Titi_2016abridged,Farhat_Lunasin_Titi_2016benard,Farhat_Lunasin_Titi_2016_Charney,Farhat_Lunasin_Titi_2017_Horizontal,Foias_Mondaini_Titi_2016,Foyash_Dzholli_Kravchenko_Titi_2014,GarciaArchilla_Novo_Titi_2018,Gesho_Olson_Titi_2015,GlattHoltz_Kukavica_Vicol_2014,Ibdah_Mondaini_Titi_2018uniform,Jolly_Martinez_Olson_Titi_2018_blurred_SQG,Jolly_Martinez_Titi_2017,Jolly_Sadigov_Titi_2015,Larios_Pei_2017_KSE_DA_NL,Lunasin_Titi_2015,Markowich_Titi_Trabelsi_2016,Mondaini_Titi_2018_SIAM_NA,Pei_2018,Rebholz_Zerfas_2018_alg_nudge}.  
Computational experiments on the AOT algorithm and its variants were carried out in the cases of the 2D Navier-Stokes equations \cite{Gesho_Olson_Titi_2015}, the 2D B\'enard convection equations \cite{Altaf_Titi_Knio_Zhao_Mc_Cabe_Hoteit_2015}, and the 1D Kuramoto-Sivashinsky equations \cite{Lunasin_Titi_2015,Larios_Pei_2017_KSE_DA_NL}. 

\subsection*{AOT Data Assimilation Algorithm}
Here, we describe the general idea of the AOT algorithm. Consider a given dynamical system:
\begin{align}
\frac{d}{dt}u = F(u,t), \qquad u(0) = u_0, \label{ref}
\end{align}
where $F$ is a possibly non-linear, possibly non-local differential operator. In applications, we require the system to be globally well posed, and it is typically assumed to have a finite-dimensional global solution. The AOT data assimilation algorithm is given by:
\begin{align*}
\frac{d}{dt}v = F(v,t) + \mu(I_h(u)-I_h(v)),\qquad v(0) = v_0,
\end{align*}

Here, $\mu > 0$ is a constant relaxation parameter, and $I_h(f)=I_h(f;X)$ denotes the piece-wise linear interpolation of $f$ at gridpoints $X=\{x_0,\ldots x_k\}$ for a grid with minimum length scale $h$, and $y$ is the solution to \eqref{ref}. In Section \ref{secSweeping}, we will consider the case when $X$ moves in time at uniform speed, and the points are closely clustered in space, mimicking a moving measurement device, such as a probe.

The number of data assimilation nodes (i.e. grid points) associated with $h$ is an important parameter for data assimilation. If the grid has too few points, the solution may not converge, but having more points than required increases computational cost. In addition to this, in real world scenarios (such as weather prediction) where these grid points are marked by sensors, such as weather monitoring devices, minimizing the number of grid points reduces the financial cost of sensor production and placement. This study examines whether one can achieve the same level of convergence using fewer points by using time-dependent data assimilation nodes. Physically this may be interpreted as moving a probe in an experiment, or mounting sensors on a moving vehicle, aircraft, satellite, etc.

The main finding of this study is the potential for using a non-stationary grid of measurement points for data assimilation. It was found that, for $\mu$ large enough (for which the simulation is stable), that in parameter ranges corresponding to more dynamically active solutions (i.e., regimes of small diffusion) far fewer nodes are required for data assimilation via a cluster of data assimilation points that move in time, which we refer to as a ``sweeping probe'' below.

\subsection*{The Allen-Cahn Equation}
The specific equation used in this study is the 1D Allen-Cahn equation (also referred to as the Chafee-Infante equation) on a bounded spatial domain (0,1):
\begin{align}\label{AC}
u_t - \nu u_{xx} &= u - \alpha u^3,\\
u(x,0) &= u_0(x),
\end{align}
with homogeneous Dirichlet boundary conditions, i.e. $u(0,t) = u(1,t) = 0$.  Here $\alpha>0, \nu>0$ are physical parameters.  In the present work, we consider the AOT algorithm adapted to this equation.  Namely, we consider the equation
\begin{align}
v_t - \nu v_{xx} &= v - \alpha v^3 + \mu(I_h(u)-I_h(v)),\\
v(x,0) &= v_0(x),
\end{align}
where $u$ is a solution to \eqref{AC}, $\mu>0$, and $I_h$ is a linear interpolation operator at nodal grid points with spacing $h$.  In the present work, we take $v_0\equiv0$ for all our simulations.  A detailed analytical study of the algorithm, and the extension of the analysis to non-static (i.e, moving) grids proposed here, will be the subject of a forthcoming work.  

\begin{figure}[htp!]
\centering
\includegraphics[scale=0.15]{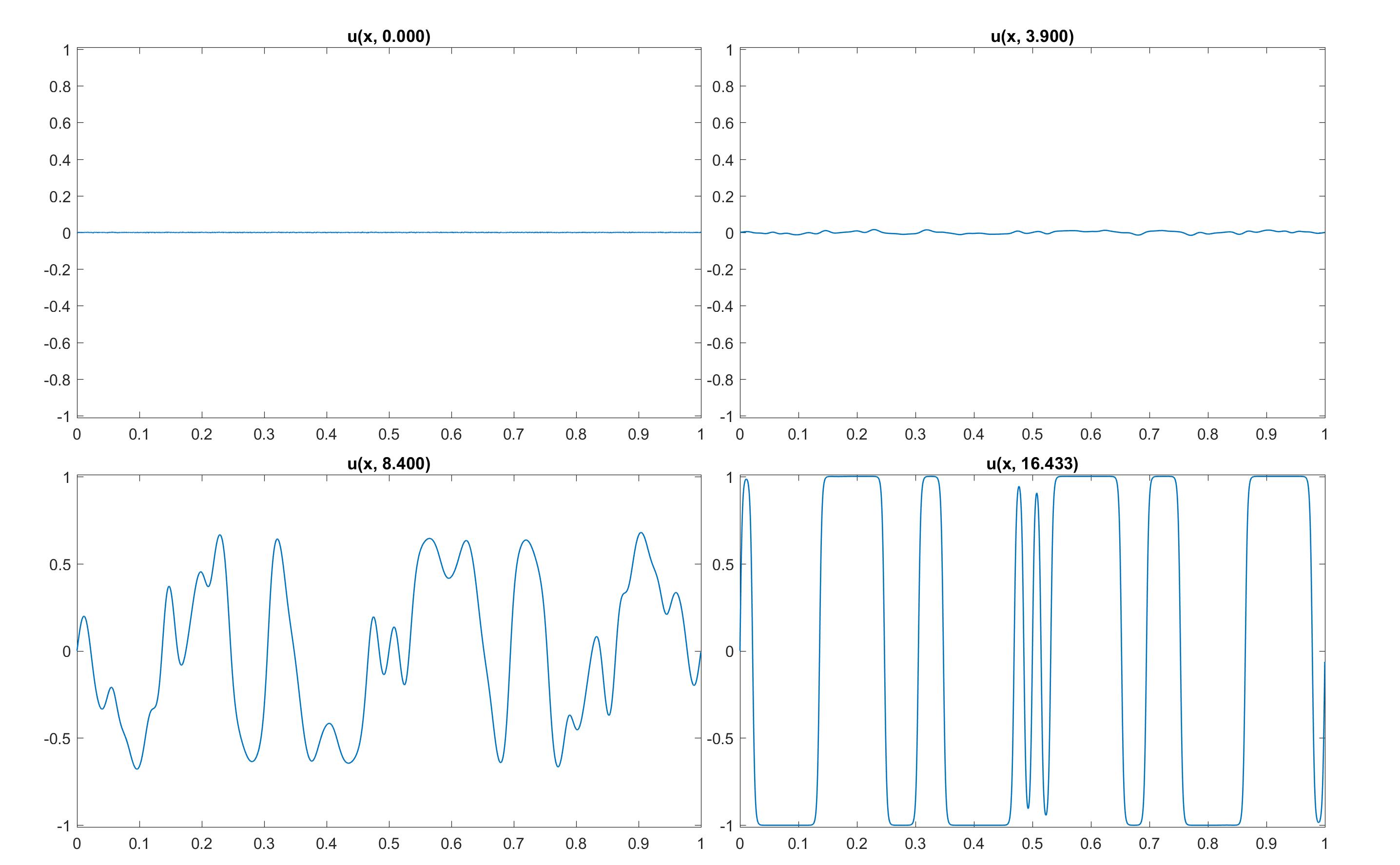}
\caption{Typical development of Allen-Cahn Equation computed solution over time for $\alpha = 1$, $\nu={7.5e-6}$. Initial small perturbations from zero (Top Left). Initial growth of small perturbations (Top Right). Development of transition layers and structures (Bottom Left). Transition layers and structures fully developed (Bottom Right).}
\label{fig:TypicalDevelopment}
\end{figure}

The Allen-Cahn equation is a chaotic system which has been studied extensively in the literature (see, e.g. \cite{Cohen_2000,Ward_1996,Shen_2010,Grant_1995} and the references therein). The assertions below can be found in, e.g. \cite{Ward_1996,deMottoni_1990,Chen_1992}. Solutions to this equation have multiple phases in their evolution. 
Typically, initial data is given as a small perturbation about zero. As the solution evolves, it goes through an initial inflationary period, marked by the development of metastable structures and transitions layers. Once metastable structures reach amplitudes of $\frac{1}{\sqrt{\alpha}}$, the solution is in a quasi-stable state with small transition layers in between structures (see, e.g., Figure \ref{fig:TypicalDevelopment}). (Note: We use exponential notation in much of the paper.  For example, we write $7.5e-6$ to mean $7.5\times 10^{-6}$.)
After the initial inflationary period, $\alpha$ determines the amplitude of the solutions, while $\nu$ influences the number of structures in the solution \cite{Robinson_2001}. The evolution of the solution at this point is metastable, in the sense that the solution will be nearly motionless for a period of time, and then very rapidly, a structure (typically the structure of smallest length scale) will disappear, being absorbed into the larger structures adjacent to it.  
This process decreases the minimal length scale of the structures, until eventually the structures stabilize, and the solution is close to the global attractor of the system.  See, e.g., \cite{Robinson_2001}, Section 11.5, for a discussion of the attractor and the large-time dynamics of the system.  
For the purposes of this study, we are most interested in the second stage, where transition layers develop, as its dynamics have certain statistically reproducible features, but each particular trajectory (in particular, the location, size, and number of structures within the trajectory) is highly sensitive to small perturbations of the initial data and the parameters of the system. However, certain features, such as the minimum length scales of the structures that appear shortly after the initial development phase, depend only on the parameters of the system (i.e., $\alpha$, $\nu$, and the length of the domain) rather than the initial conditions (see, e.g., \cite{Robinson_2001}).

\subsection*{Static vs. Non-Static Measurement Points}

One issue in working with static measurement point locations is that it may be difficult to know where to optimally place the measurement devices (which we call \textit{sensors} below).  Another difficulty is in knowing the total amount of sensors required to achieve convergence. When modeling the Allen-Cahn equation we found that in order to achieve convergence to the reference solution using AOT data assimilation the placement of the sensors was vital. We observed in computational experiments (described below) that if sensors were not placed in certain critical regions (corresponding to wave structures and transition layers between them), solutions did not converge to the reference solution. Unfortunately, the locations of the regions are determined by the dynamics of the solution. Moreover these locations vary chaotically in the sense that they are highly sensitive to perturbations in the initial data and the parameters for the system. Thus, \textit{a priori} determination of the critical sensor locations appears to be impractical. Obtaining convergence using a static grid therefore appears to necessitate a uniform covering of the entire domain with a grid of sensors fine enough to cover all possible locations where critical regions might arise.  To get around this difficulty, in the present work, we investigate the use of time-dependent sensor grids. We replace the static grid with a sweeping probe (a moving cluster of data assimilation points).  This probe moves across the domain as time evolves, thereby entering every potential location where a critical region might develop.
We found that using this strategy, we could achieve the same rate of convergence using fewer nodes.  In fact, this improvement was most dramatic in parameter ranges corresponding to more chaotic solutions (i.e., regimes of small diffusion), differing by as much as an order of magnitude in the smallest viscosity regimes we tested.

This technique may have applications to other settings, such as turbulent flow.  We will investigate such applications in a forthcoming work.

\subsection*{Applications to Parameter Estimation}

By examining the relationship between convergence of various static grids, we find an estimate for a minimum length scale $\lambda$. Experimentally, we discovered that the convergence of the static grid was based on the placement of the sensors in critical regions. This observation led us to a heuristic argument, presented in Section \ref{secInv}, that the minimum number of sensors required for convergence,$m_h$, is directly related to $\lambda$. We then experimentally find a relationship between $\nu$ and $m_h$ in the case of a uniform static grid.  Combining all of this information algebraically gives an estimate for $\lambda$ directly in terms of $L$ and $\nu$. Such an approach may be a useful for solving certain inverse problems.
%

\subsection*{Organization of the paper}
The paper is organized as follows.
In Section \ref{secNum} we describe the numerical methods we use to simulate the PDE, and also the methodology of our computational investigation. 
For the purposes of comparison with the standard AOT algorithm on a uniform static grid, in Section \ref{secUnifGrid} we run several simulations in this setting. 
In Section \ref{secSweeping} we present our main results on the case of time-dependent measurement points.  We show that, in the simulation regimes we considered, the number of grid points needed for convergence using a sweeping probe is far smaller than the number needed for convergence of the algorithm using a uniform grid over a wide range of parameters.
In Section \ref{secInv} we provide an application of the AOT algorithm to an inverse problem.  Namely, we show that the minimum length scale of structures in solutions to the Allen-Cahn equations that appear shortly after the initial transient phase.

\section{Numerical Methods}\label{secNum}
\subsection*{Solving the PDEs}
We solve this equation on a uniform grid of $N$ points distributed uniformly on $\left[0,1\right]$ using uniform discrete time-steps which are stable for the implicit/explicit scheme we used, according to \cite{Eyre_1997}. To ensure that initial data had a fully-resolved Fourier spectrum we used initial data in the form:
\begin{align}\label{initial_data}
u_0(x) = \sum_{k=1}^{\frac{N}{4}}a_k\sin(2\pi kx),
\end{align}
with $a_k$ determined randomly using a normalized Gaussian distribution, and $u_0$ then rescaled by a constant so that $\|u_0\|_{L^2}=1e-2$. This solution can be thought of as periodic and odd on $[-1,1]$ or as satisfying homogeneous Dirichlet boundary conditions on $[0,1]$, since the equation preserves the periodicity and oddness of the initial data.

This study uses a semi-implicit convex splitting scheme from \cite{Eyre_1997,Eyre_1998} to solve for the reference solution and the simulated data assimilation solution at every time-step. The numerical scheme we use was derived by \cite{Eyre_1997,Eyre_1998} as a stable implicit/explicit scheme:
\begin{align*}
{U_k}^{n+1} - {U_k}^n = dt(\nu D_{xx} + 1+2\alpha){U_k}^{n+1} + dt(-2\alpha - \alpha ({U_k}^{n})^{^2}){{U_k}^n},
\end{align*} 
where ${U_k}^{n} = u(x_k,t_n)$ where $dx = \frac{1}{N}$, $x_k = k dx$, $k = 0,1,2,...,N-1$, $dt = 1e-3$, $t_n = n dt$. 
$D_{xx}$ is a centered-difference approximation of the operator $\frac{\partial^2}{\partial x^2} $, here $\frac{\partial^2 {U_k}^{n+1}}{\partial x^2}$ is given by the second-order finite difference approximation $\frac{{U_{k-1}}^{n+1} - 2{U_{k}}^{n+1} + {U_{k+1}}^{n+1}}{dx^2} $. This results in the below system which needs to be solved at every time-step:
\begin{align}\label{eq:Num_Scheme}
(1 - dt(\nu D_{xx} + 1+2\alpha)){U_k}^{n+1} =  (1+dt(-2\alpha - \alpha( {U_k}^{n})^{^2})){{U_k}^n}.
\end{align}
This method, being semi-implicit, is more stable than fully explicit methods in the sense that it allows for a much larger time-step than fully explicit methods \cite{Eyre_1997}. In addition to this, the resulting matrix is tri-diagonal, so solving this system has a time complexity of only $O(N)$ using the Thomas algorithm.  In every simulation in this work, we use time-step $\Delta t = 1e-3$ for stability, and spatial resolution $N = 2^{12}=4096$ (i.e., $\Delta x = 2^{-12}\approx 2.441e-4$), which are sufficient to resolve the spatial scales in our simulations (i.e., below the cubic aliasing cut-off number $N/4$) to machine precision over the range of $\nu$-values and $\mu-values$ we consider.  Note that since we handle the data assimilation term explicitly, a CFL condition arises, requiring $\Delta t\leq 2/\mu$, which is satisfied in all of our simulations.

In order to apply data assimilation, we solve the following system at each time-step:
\begin{align}\label{eq:V_num}
&
(1 - dt(\nu D_{xx} + 1+2\alpha)){V_k}^{n+1} 
\\\notag&\quad  =
({V_k}^n+dt((-2\alpha{V_k}^n - \alpha ({V_k}^{n})^{^3}) +\mu(I_h(U^n) - I_h(V^n))),
\end{align}
where we treat the data assimilation term explicitly, and $I_h(U^n)$ is the result of algorithm \eqref{eq:Num_Scheme}. We initialize this system with ${V_k}^0 = 0$ for all $k$.

It is important to note that the scheme above generates a numerical approximation for $u$ and $v$. These numerical approximations will be denoted as $\widetilde{u}$ and $\widetilde{v}$ respectively. One must also note that the term $U^n$ in \eqref{eq:V_num} is referring to the numerical approximation $\widetilde{u}$ evaluated at time $t = t_n$, i.e. $U^n = \widetilde{u}(t_n)$.

The data assimilation scheme tested in this study used the standard static uniform grid for the interpolation operator $I_h$. Here, $I_h(f) = I_h(f;X)$ is the piecewise linear interpolation of $f$ based on the gridpoints $X$, associated with length-scale $h$. For the uniform grid case, we take $h$ to be an integer multiple of $\Delta x$.  For the moving probe case, $h$ does not have a particular meaning, but we retain the notation $I_h$ for consistency.

\section{Uniform Static Grid}\label{secUnifGrid}

In order to understand the performance of the standard AOT algorithm in the context of the Allen-Cahn equations, and also for comparison with the moving probe simulation that we consider in the next section, we consider in this section the algorithm in \eqref{eq:V_num} with uniformly spaced measurement points which are static in time.

To run the simulations, we initialized the reference solution by choosing random periodic initial data with small amplitude for $u_0$, chosen as in \eqref{initial_data}, and using the numerical scheme in \eqref{eq:V_num} and \eqref{eq:Num_Scheme} compute the approximate numerical solution $\widetilde{u}$.  We evolve the system until it develops metastable structures. That is, the system evolves until there is a maximum value being within 20\%  of $\frac{1}{\sqrt{\alpha}}$, the maximum amplitude of the system. This process can be seen in Figure \ref{fig:TypicalDevelopment}. We then initialize the data assimilation solution $\widetilde{v}$ with identically zero initial data and reference solution $\widetilde{u}$ using the evolved system state. We next create a uniform static grid consisting of $m_h$ uniformly distributed points with length scale $h$. The reference solution and the simulated AOT algorithm solution then run in tandem using a static uniform grid configuration with approximately $h^{-1}$ points. We allow these systems to develop to time $t_s$, where $t_s$ is the time after the initial development of metastable structures for the reference solution (we took $t_s:=\max\{t\leq 10:\norm{\widetilde{u}(t)}_L^{\infty}\leq0.8/\sqrt{\alpha}\}$). This was done because the time to develop these structures vary with the parameters $\nu$ and $\alpha$, making a uniform ending time a possibly biased comparison for convergence of data assimilation solutions. 

We take the error at time $t$ to be given by $\norm{\widetilde{u}(t)-\widetilde{v}(t)}_{2}$ is used as an approximation of the error. The second phase is repeated using a binary search algorithm to determine the minimum value of $m_h$ for which there is sufficient convergence. 
We take \textit{sufficient convergence} to mean that, by a certain time $t_*$, the $L^2$-error in the solution is below a tolerance of $5 e-14$, slightly larger than MATLAB's machine epsilon, namely $2.2204e-16$.  That is, we will say a solution is converged at time $t_*$ if 
\begin{align}\label{sufficient_convergence}
\norm{\widetilde{u}(t_*)-\widetilde{v}(t_*)}_{L^2} \leq 5e-14.
\end{align}
Here, we take $t_* = 50$.  To approximate the value of $m_h$, we used a binary search, varying the number of nodes in a uniform mesh on the interval $\left[1,N\right]$; that is, starting with a lower number of nodes equal to $1$, and an upper number of nodes equal to $N$.  We then test $\lfloor N/2\rfloor$ uniformly spaced nodes for sufficient convergence.  We repeat the process on the interval $\left[1,\lfloor N/2\rfloor\right]$ if the test is successful (i.e, sufficient convergence is attained), testing with $\lfloor N/4\rfloor$ grid points, or on the interval $\left[\lfloor N/2\rfloor,N\right]$ testing with $\lfloor 3N/4\rfloor$ grid points, if the test was unsuccessful, and so on, until the interval is length 1, and we set $m_h$ equal to the upper endpoint of the interval.

We note that placement of the nodes will change as the number of nodes changes.  For instance, if one uses three evenly-spaced nodes on the interval $[0,1]$, the nodal locations will be different that if ones uses four evenly-spaced nodes on the same interval.  We note that this could lead to inconsistent results if, e.g., a certain solution has transition layers that happen to align with these nodes, which could potentially leads to spurious convergence rates.
Therefore, we run 
10 trials based on different randomly generated initial data for each fixed value of the viscosity $\nu$.
The results of some of these trials can be seen in Figures \ref{fig:MinimunGrid500} and \ref{fig:MinimunGrid1000}. Some variance in $m_h$ was found by altering the value of $\mu$, with large $\mu$ decreasing $m_h$.

\begin{figure}
	\centering
	\includegraphics[scale=0.15]{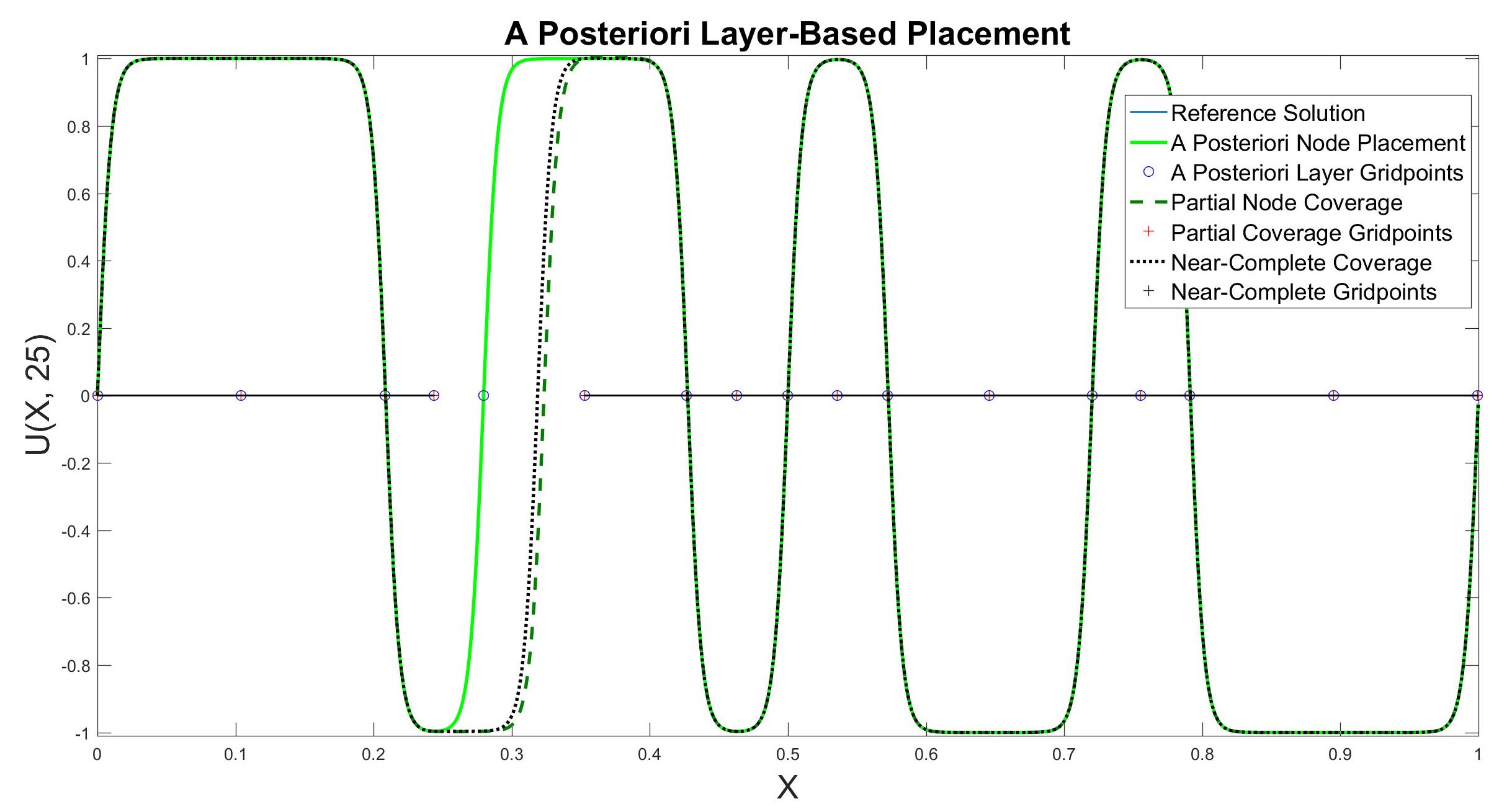}
	\caption{Data assimilation can converge or diverge based on an interval not being covered $\alpha = 1$, $\nu={5e-5}$, at $t_s=25.0$.  Note that near $x = 0.3$, the simulation is not converged.}
	\label{fig:Placement}
\end{figure}

\begin{figure}
	\centering
	\includegraphics[scale=0.15]{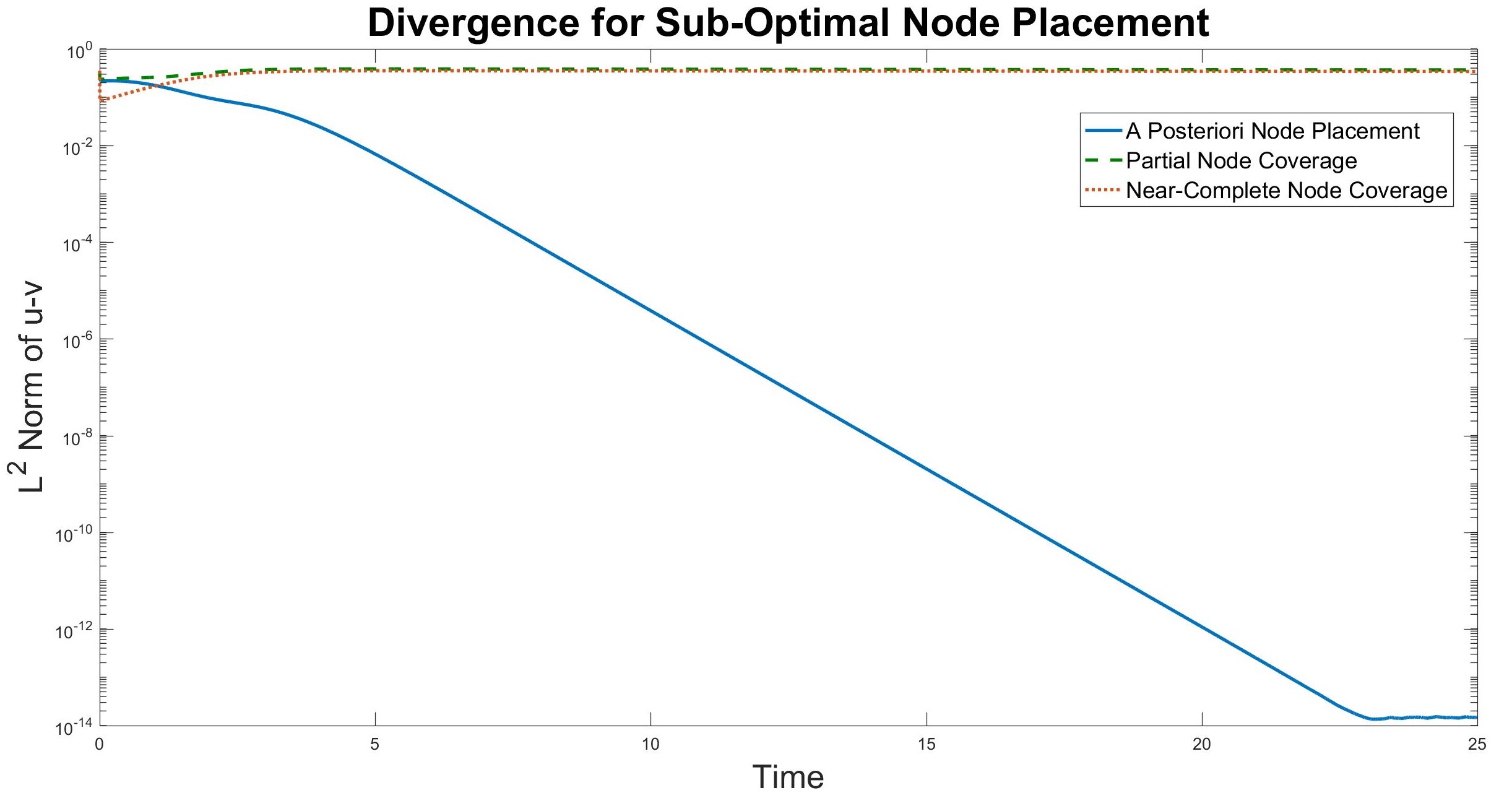}
	\caption{Error associated with Figure \ref{fig:Placement} vs. time  (log-linear plot).}
	\label{fig:PlacementError}
\end{figure}

To get a more optimal estimate for comparison, we worked backwards to generate $m_h$ values using \textit{a posteriori} knowledge of the location of the transition layers. While not useful in practice, due to the need for \textit{a posteriori} information, using an \textit{a posteriori} layer-based placement in some sense gives an ``optimal'' or ``minimal'' number of measurement nodes. That is, it might be thought of as a ``best case'' scenario for static grids which are not necessarily uniform.  To investigate this we conducted three separate experiments, the results of which can all be seen in Figure \ref{fig:Placement} and Figure \ref{fig:PlacementError}.  We then performed data assimilation with a static grid configured using an \textit{a posteriori} layer-based placement strategy, namely, one data assimilation node was placed in each transition layer, structure, and at each endpoint.  Then using the same initial data we repeated this experiment, but with one data assimilation node in one transition layer removed. Finally we repeated the previous experiment but with data assimilation nodes at every other point on the domain except in one interval containing a transition layer. 

The estimated error found in these studies can be seen in Figure \ref{fig:PlacementError}. As one can see the only trial that resulted in convergence to the reference solution was the \textit{a posteriori} layer-based placement strategy. It appears there is a visible difference in the end-state of each solution as well (see Figure \ref{fig:Placement}). Similar results were found when excluding data assimilation nodes from a single structure as well. Based on these experiments it appears that, for this initial data, that the positioning of the data assimilation nodes is an important factor for the eventual convergence of the system to the reference solution. 



\section{Sweeping Probe Data Assimilation}\label{secSweeping}
\noindent
The major aim of this section is to investigate the use of a sweeping probe (i.e. a cluster of data assimilation points which move in time, eventually covering the entire domain) for the AOT data assimilation algorithm. That is, we consider the following data-assimilation scheme on the 1D domain $[0,L]$:
\begin{subequations}
 \begin{align}
v_t - \nu v_{xx} &= v - \alpha v^3 + \mu(I_h(u;X_t)-I_h(v;X_t)),\\
X_t&:= \{x_1 + ct,x_2 + ct,\ldots,x_{m_h} + ct\} (\text{mod} L)\\
v(0,t) &= v(L,t) = 0,\\
v(x,0) &= v_0(x),
\end{align}
\end{subequations}
where $x_j = jh$ for some $h>0$ and $j=1,2,\ldots,m_h$, and $c>0$ is a constant.  As above, we take $L=1$.  The ``(mod L)'' notation means at the endpoints of the domain, the probe wraps around periodically.  We do this to avoid oversampling in time near the boundaries.  We only consider constant speed $c$ here, but in a future work, we will study non-uniform movement of the probe, which may also depend on the input from local measurements.  

\begin{figure}
	\centering
	\includegraphics[scale=0.15]{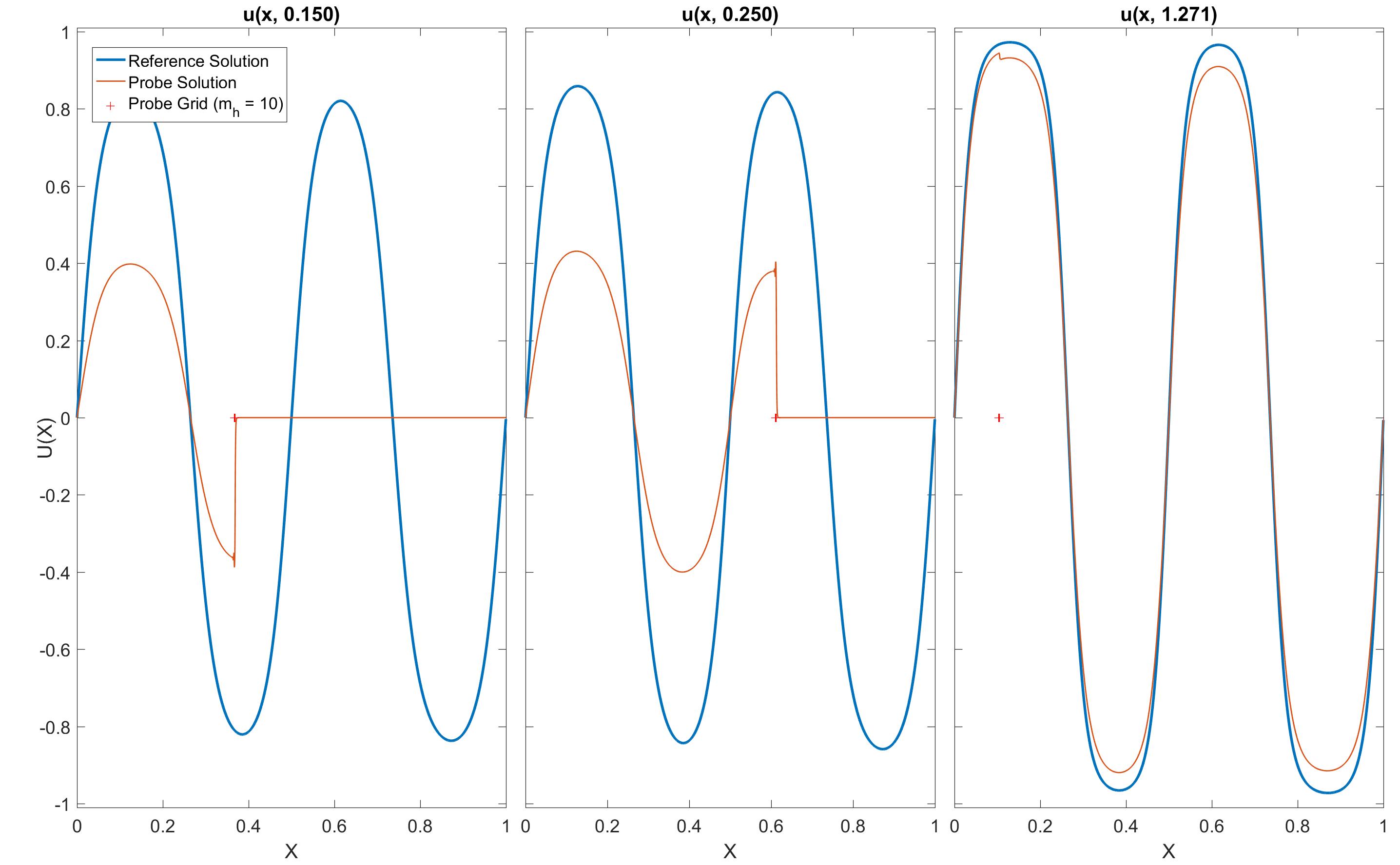}
	\caption{
	Snapshots of the sweeping-probe simulation in time, with data points clustered around a moving point located at the red "+". The assimilation simulation (orange curve) grows toward the true solution (blue curve) in time.   Parameters: $\alpha = 1$, $\nu={5e-4}$, $\mu = 500$, $c = 10 \frac{\Delta x}{\Delta t}$ and $m_h = 10$ at times $t = .15$, $t=.25$, and $t=1.271$.
}
	\label{fig:ProbeDev}
\end{figure}

We found in our simulations that parameter regimes corresponding to sufficiently small values of $\nu >0$ require far fewer nodes for data assimilation via a sweeping probe than for a uniform placement of nodes (see Figures \ref{fig:MinimunGrid500} and \ref{fig:MinimunGrid1000} below). To investigate this, we repeated the process used for finding the minimum number of data assimilation nodes for a uniform grid (the binary search algorithm outlined in Section \ref{secUnifGrid}), but instead of using a static uniform grid, we use a small cluster of $m_h$ consecutive points spaced apart by the finest length scale $h=\Delta x$.  In our simulations, $m_h$ was typically so small that the total length of the probe was $|x_{m_h}-x_1|< 0.005L$. This cluster of points moves at a constant velocity $c$. At each time-step we shift every point to the right by $c\cdot\Delta t$ units, where $c$ is chosen so that $X_t$ aligns with the underlying spatial mesh at every discrete time-step.   Figure \ref{fig:ProbeDev} shows the typical development of the probe solution over time.  The data assimilation solution is initially to zero at time $t=0$.  At time $t=0.150$, the probe has moved from $x=0$ to $x=0.4$.  At time $t=0.250$, the probe has moved further along in space, and has continued to feed data into the simulation.  At time $t=1.271$, the probe has wrapped periodically around the domain twice, and is continuing to feed data into the simulation, which has nearly converged to the true solution.  For larger times (not pictured), the  assimilated solution will converge to the true solution up to machine precision.

The probe moves through the domain, sampling every point. As the probe passes through the domain it produces nearly discontinuous solutions, strongly forcing the simulation towards the reference solution as it passes through the domain. After the first pass through domain, the structure of the sweeping probe solution mimics the basic structure of the reference solution, i.e. the structures and transition layers form in the correct locations with smaller amplitude than in the reference solution. With each pass through the domain the probe solution more closely approximates the reference solution, eventually converging to within machine precision for sufficiently large $m_h$ values. It is important to note that the sweeping probe uses the finest length scale on the discretized domain. 
In this context, the number of nodes in the probe is proportional to the length of the probe.


To approximate the value of $m_h$ for the moving probe case, we used a binary search, as in the uniform static grid case.  Also, for the sake of consistency with that case, we also ran 10 trials based on different randomly generated initial data for each fixed value of the viscosity $\nu$. Unlike the uniform grid case, it was never observed in our tests that increasing the number of nodes in the probe negatively affected convergence rate in the sweeping probe case. 
The results of some of these trials can be seen in Figures \ref{fig:MinimunGrid500} and  \ref{fig:MinimunGrid1000}.
\begin{figure}[htp!]
	\centering
	\includegraphics[scale=0.15]{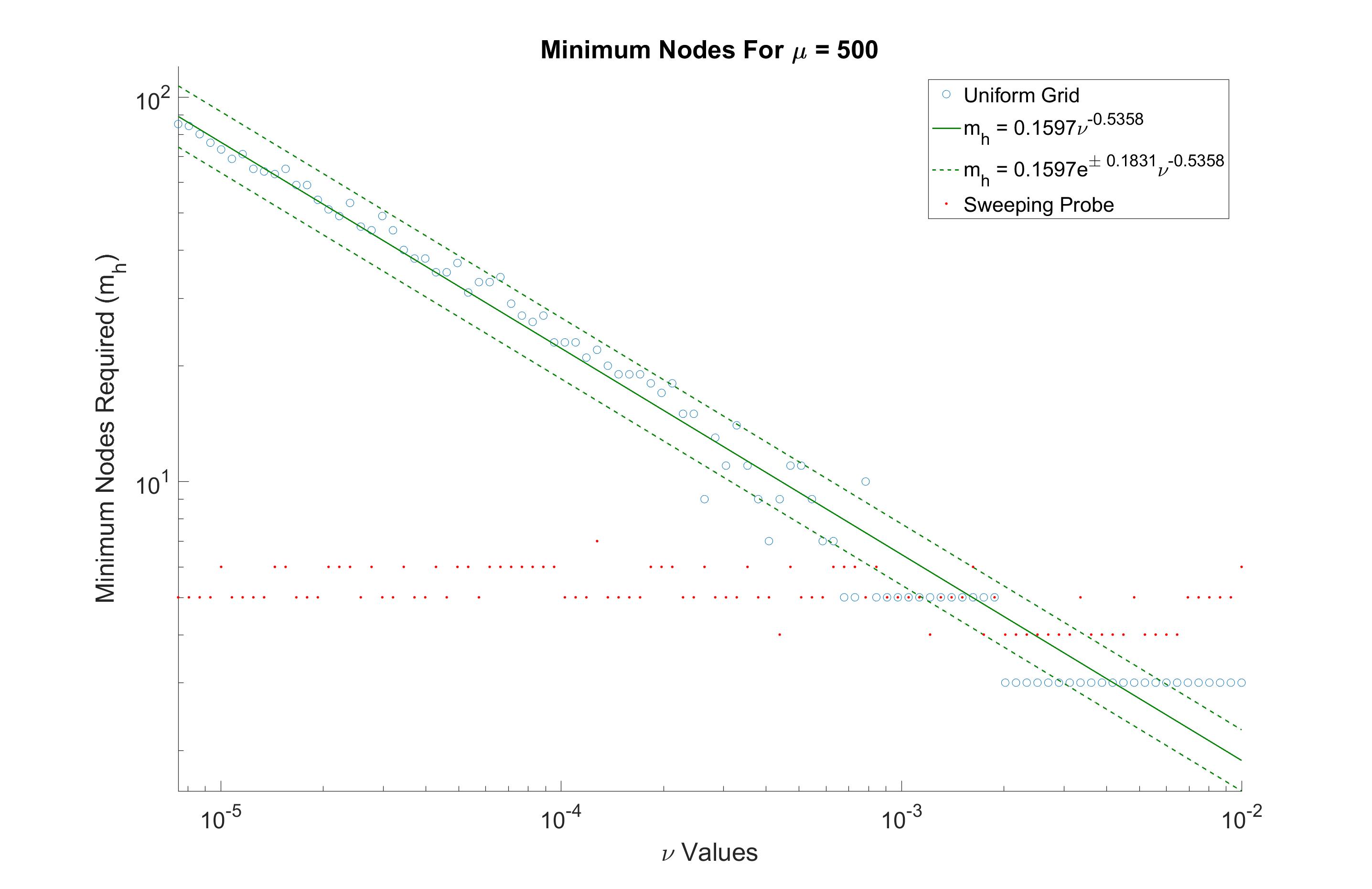}
	\caption{Minimum  number of nodes for uniform grid and sweeping probe required for convergence of $5e-15$, $t_s=50$ units after initialization with an exponential approximation for the uniform grid given by $m_h = 0.1597 e^{\pm 0.1831}\nu^{-0.5358}$.  Here, $\alpha = 1$, $\mu = 500$, and $c = 30\frac{\Delta x}{\Delta t}$. }
	\label{fig:MinimunGrid500}
\end{figure}
\begin{figure}[htp!]
	\centering
	\includegraphics[scale=0.15]{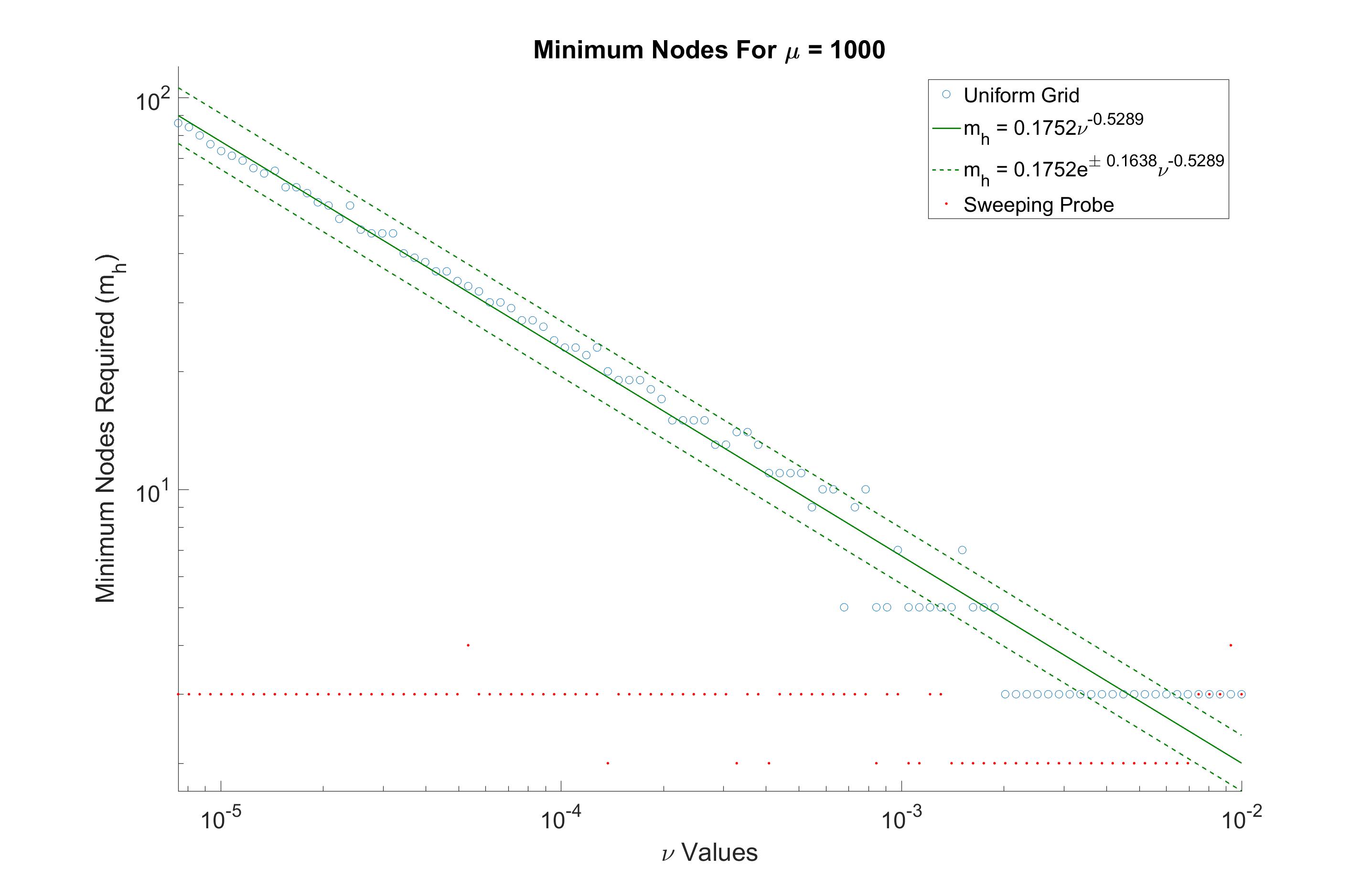}
	\caption{Minimum number of nodes for uniform grid and sweeping probe required for convergence of {5e-15}, $t_s=50$ units after initialization with an exponential approximation for the uniform grid given by $m_h =  0.1752\cdot e^{ \pm 0.1638}\nu^{ -0.5289}$.  Here, $\alpha = 1$, $\mu = 1000$, and $c = 30\frac{\Delta x}{\Delta t}$. }
	\label{fig:MinimunGrid1000}
\end{figure}

We found that for a sweeping probe with a constant velocity $c = 30\frac{\Delta x}{\Delta t}$ (i.e. the probe moves 30 gridpoints every time-step), and with $\mu = 1000$, that the sweeping probe needed fewer data assimilation nodes for all tested values of the viscosity $7.5e-6\leq \nu\leq 7e-3$.  For $\nu> 7e-3$, the sweeping probe needed the same number of data assimilation nodes as the uniform grid in four of the five trials.  In only one of the 100 viscosity values tested did the sweeping probe need more nodes than the uniform grid, and in this case, it only needed one more node.
(see Figure \ref{fig:MinimunGrid1000}). Similar results were found for smaller values of the velocity  $c$ as well.

\begin{figure}
	\centering
	\includegraphics[scale = .15 ]{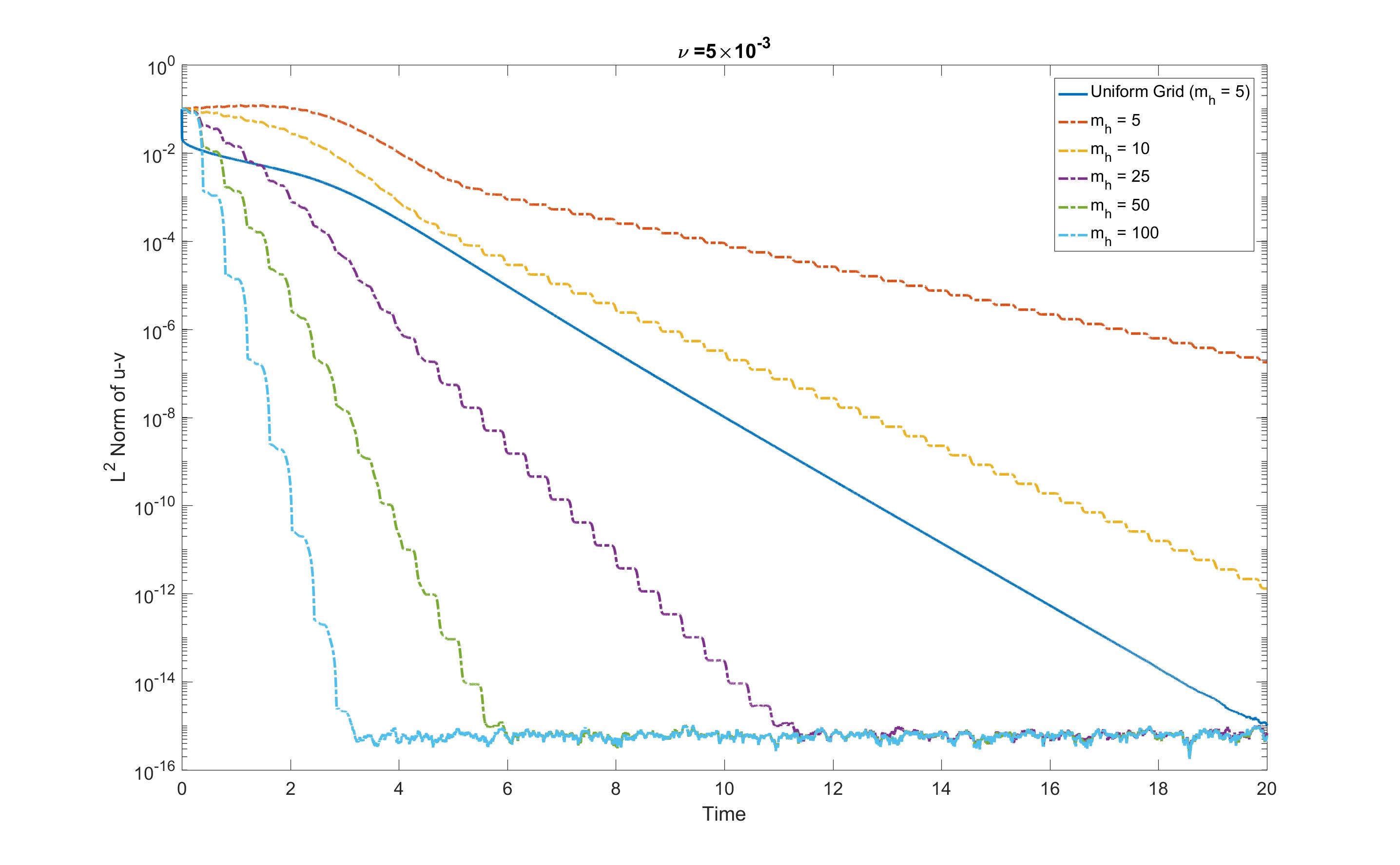}
	\caption{Comparison of error over time for static uniform grid and sweeping probe of various lengths (log-linear plot).  Here, 
	$\nu = 5e-3$, $\alpha = 1$, $\mu = 500$, and $c = 10 \frac{\Delta x}{\Delta t}$.}
	\label{fig:CarGridComp2}
\end{figure}

\begin{figure}
	\centering
	\includegraphics[scale = .15 ]{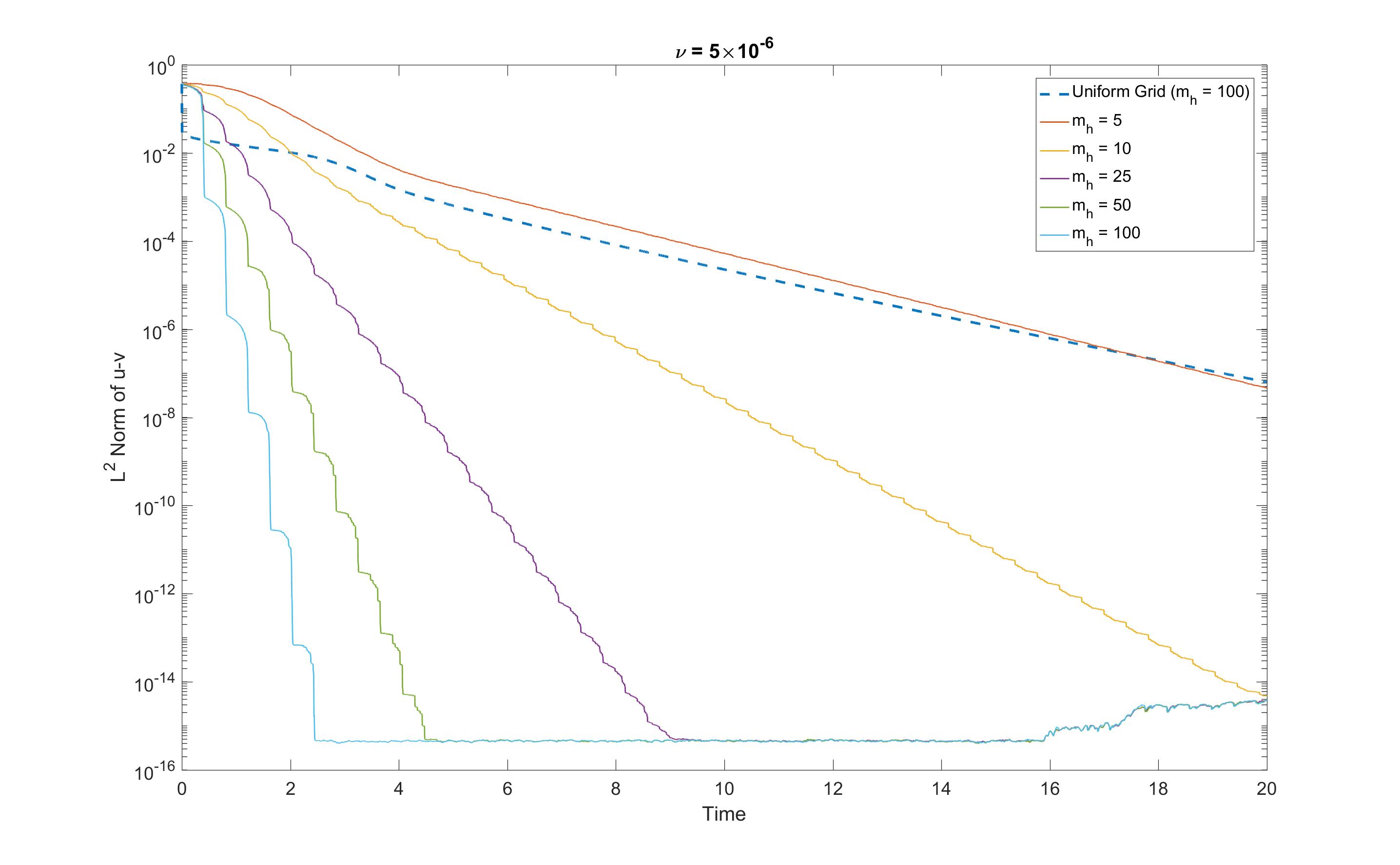}
	\caption{Comparison of error over time for static uniform grid and sweeping probe of various lengths (log-linear plot). Here,
		$\nu = 5e-6$, $\alpha = 1$, $\mu = 500$, and $c = 10 \frac{\Delta x}{\Delta t}$.}
	\label{fig:CarGridComp1}
\end{figure}

Another important metric to examine when comparing performance of the uniform static grid case to the sweeping probe case is the convergence rate of the error over time. To investigate this, we ran multiple simulations for probes of varying lengths and tracked the $L^2$-error over time. We then compared these convergence rates to those of a uniform static grid initialized with the same initial data and a sufficiently large number of data assimilation nodes required for convergence.  
The results of this can be seen for $\nu = 5e-6$ and $\nu = 5e-3 $ in Figure \ref{fig:CarGridComp2} and Figure \ref{fig:CarGridComp1}, respectively. 

For large viscosity (e.g., $\nu = 5e-3$ in Figure \ref{fig:CarGridComp2}), the uniform grid case outperformed the moving probe case in the sense that the moving probe needed more nodes to converge to machine precision at the same time.  However, in the case of smaller viscosity (e.g., $\nu = 5e-6$ in Figure \ref{fig:CarGridComp1}), the situation is reversed, and we see that the moving probe case outperforms the uniform grid case.  Indeed, its convergence rate, even with 10 nodes is far greater than the convergence rate of the uniform grid case with 100 nodes--an order of magnitude improvement in the required number of nodes.  In addition, this rate improves rapidly as the number of nodes increases toward 100.


\begin{remark}
 The error curves for the sweeping probe trials in Figures \ref{fig:CarGridComp2} and \ref{fig:CarGridComp1} follow a descending stair-step pattern. We note that the length of each step corresponds to the time required for the probe to pass through the entire domain.
\end{remark}

\begin{figure}
  	\centering
 	\includegraphics[scale = .15 ]{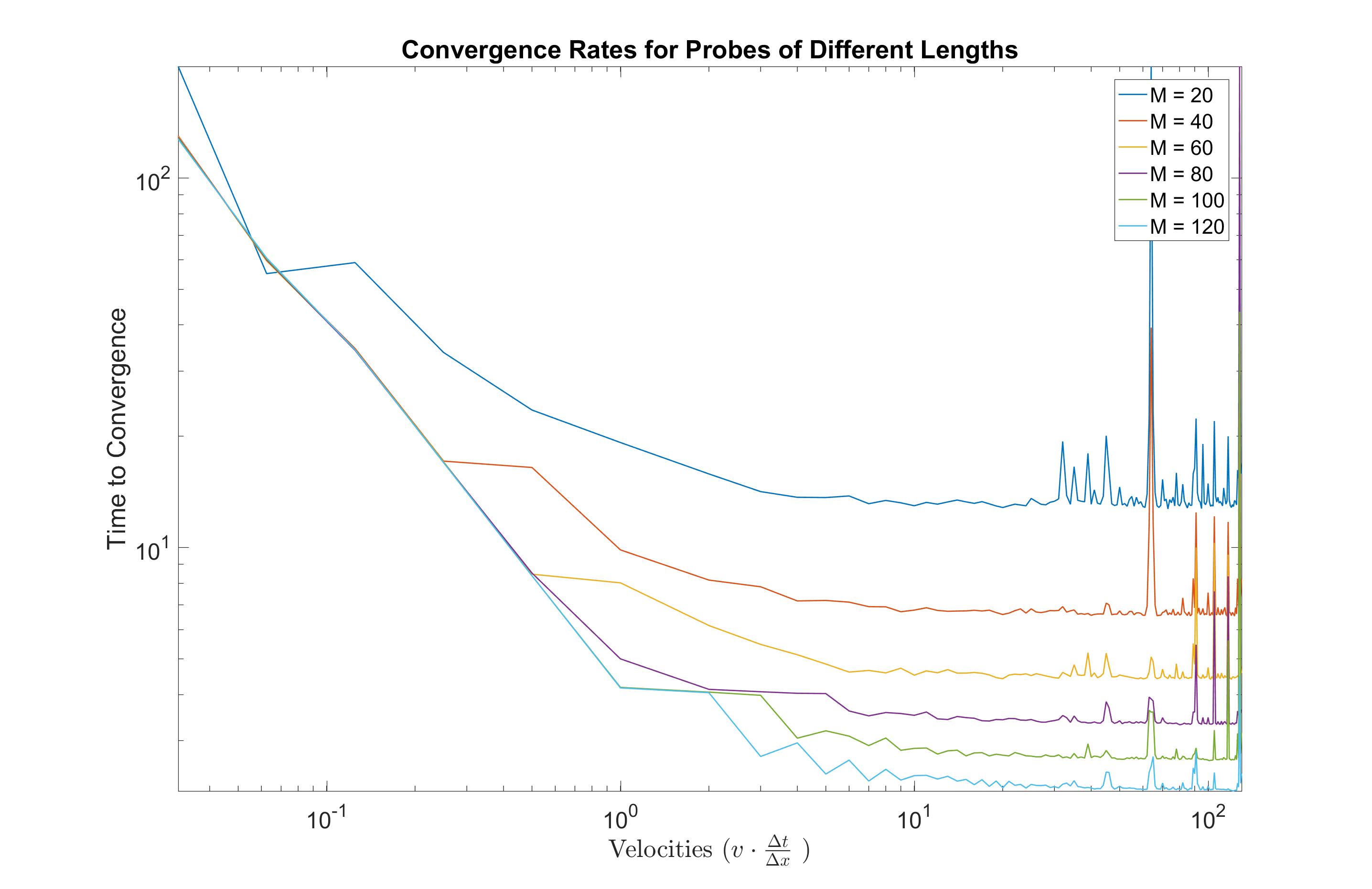}
 	\caption{Times for convergence to $5e-15$ with probes of differing velocities and sizes (log-log plot).  Here, $\nu = 7.5e-6$, $\alpha = 1$, and $\mu = 300$.}
 	\label{fig:VelocityTimes}
 \end{figure}

We also investigated the effect of the velocity of the probe on time to convergence. The velocities in question are static velocities related to the time-step, i.e. the probe moves a constant number of gridpoints each time-step. We also set up another experiment where we initialized a sweeping probe of a variable number of sensors. 
We allowed the system to develop and then tracked the amount of time it took for the sweeping probe to converge to the reference solution for a given velocity. This was repeated over 130 times using the same initial data and number of probe sensors ($m_h$), but with several different values of the velocity $c \in \left(0, 130\frac{\Delta x}{\Delta t}\right]$. The results of these experiments can be seen in Figure \ref{fig:VelocityTimes}.

In these experiments, we found that if the probes converged to the reference solution, they tended to converge at a time which  was close to the mean convergence time for all velocities, with slight perturbations from this mean. (Here we take sufficient convergence to be when $\norm{\widetilde{u}(t_*)-\widetilde{v}(t_*)}_{L^2} \leq {e-10}$ at some time $t_*>0$.)  That is to say, in general, the speed of the probe appear not have a strong effect on the convergence time.  However, we did observe two types of ``suboptimal'' velocities that occur, which we call ``frequency-locked velocities''. By suboptimal, we mean that the time to convergence was greater than 1.5 times the mean convergence time (averaged over trials), or that convergence was never achieved. 
The first type of frequency-locked velocity is $c = 0$, namely, when the probe does not move.  Here, the simulation did not converge.  Moreover, for velocities smaller than $c = 1\frac{\Delta x}{\Delta t}$, convergence was slow, and the convergence time appeared to increase without bound as $c$ tended towards zero.
The second type of frequency-locked velocity we observed was $c = a\cdot K\frac{\Delta x}{\Delta t}$, for every $a=1,2,3\ldots$. $K>0$ here is an integer dependent on $N$, the minimum length scale $\lambda$ (see Section \ref{secInv}) and the number of nodes in the probe, $m_h$. For $N = 2^{12}$ and $\nu = 7.5 e-6$, the frequency-locked velocities were observed to correspond to the value $K = 64$. It was found that the probe took significantly longer to converge or it did not converge at all at these velocities. This is to be expected due to the periodicity of the probe. Since the probe sweeps periodically around the domain of $N = 2^{12}$ points, when $c = a\cdot K \frac{\Delta x}{\Delta t}$ the probe does not sample every point in the domain, but continuously hits a specific subset of gridpoints. At these velocities, the gaps between these subsets are too large (e.g., large enough that a transition layer might never be sampled), and as expected, the simulation was never observed to converge.  Based on these observations, we conjecture that the following conditions imply that $K$ corresponds to a frequency-locked velocity: $K$ divides $N$ and $|K- m_h| < \lambda/\Delta x $. 

\begin{figure}
	\centering
	\includegraphics[scale = .15 ]{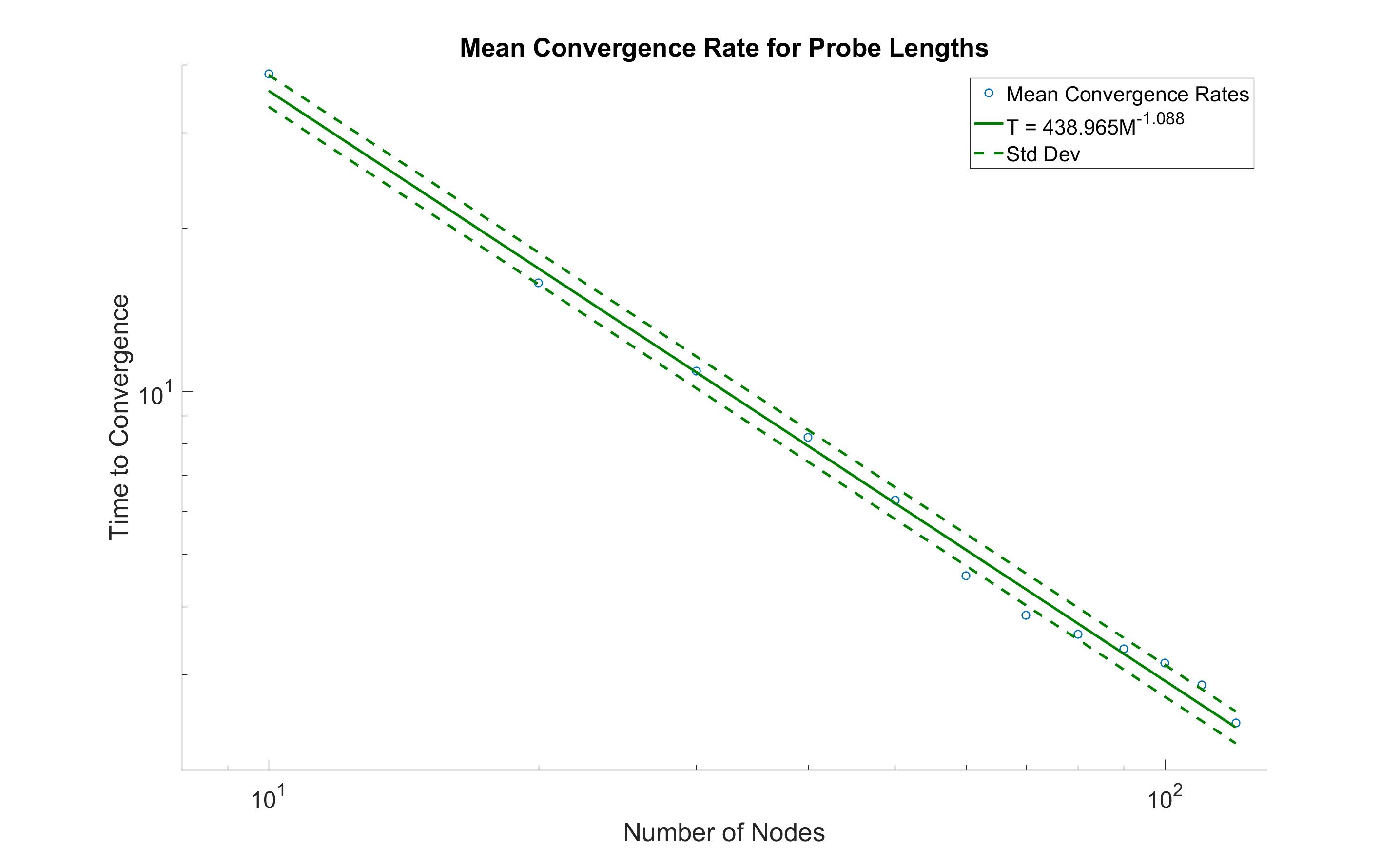}
	\caption{Mean times for probes of different lengths to converge within $5e-15$ (log-log plot), with an exponential approximation given by $T = 438.965M^{-1.088}$.  Here,
	$\nu = 7.5e-6$, $\alpha = 1$, and $\mu = 300$.}
	\label{fig:VelocityTrend}
\end{figure}

Increasing the number of nodes has the benefit of decreasing the mean convergence time for all velocities. We measured the convergence rate of the same initial data with $\mu = 300$, $\nu = 7.5e-6$, and $\alpha = 1$ and we observed that (as expected) whenever the number of nodes was increased, but the velocity remained the same, the time to convergence (to machine precision) of the simulated solution was smaller. The results of these experiments can be seen in Figure \ref{fig:VelocityTrend}. In the case of $\nu = 7.5e-6$, the mean convergence time, $T$, for a probe of size $M$, we found $T \approx 438.965M^{-1.088}$.



\section{Applications to an Inverse Problem}\label{secInv}
In the process of completing this manuscript, we noticed that some of our work in the uniform static grid case could be used to estimate the minimal length scale of solutions as they enter the second stage of their evolution, for each given diffusion coefficient.  

Repeating the experiments for the uniform grid, described in Section \ref{secUnifGrid}, at least 10 times in each case, for a range of $\nu$ values ($\nu \in [{7.5e-6},{1e-2}])$, we make an ansatz that the minimum number of data assimilation nodes required for convergence, $m_h$, and the diffusion coefficient $\nu$ are related by a power law in the form
\begin{align}\label{mh_est}
m_h \sim c_0\nu^{-p}.
\end{align}
This relationship is supported by the results of our simulations in Figure \ref{fig:MinimunGrid500}).  Moreover, using linear regression (on a log-log scale), we estimate the parameters from our data as $c_0\approx 0.1752\cdot e^{ \pm 0.1638}$ and $p\approx 0.5358$.

In all cases we found that for $m_h$ sufficiently close to $c_0\nu^{-p}$, the $L^2$-error satisfied convergence condition \eqref{sufficient_convergence}. The results of these experiments can be see in Figure \ref{fig:MinimunGrid1000}. This estimate and the corresponding error-bars were found by applying curve fitting techniques to determine an approximation of the form $m_h = a\nu^b$ (namely, linear regression on a log-log scale). Section \ref{secUnifGrid} indicated that a solution for which all the structures are smallest, and therefore uniformly distributed, should require the largest number of measurement nodes using a layer-based placement strategy. Moreover, this indicates that least number of nodes (using a layer-based placement) would be required when all the structures have the same size $\lambda$. Let us denote by $n_s$ the maximum number of structures. Of course we are only looking at the maximum number of structures that develop initially, see \cite{Robinson_2001} for the number of structure that develop as $t\rightarrow \infty$. In this case we should have $\lambda = \frac{L}{n_s} $, where $L$ is the length of the domain. Noting that $m_h = 2n_s$ or equivalently $n_s = \frac{m_h-1}{2}$, since for every structure, there will be a transition layer to the right and each structure shares a transition layer with the neighboring structure except at the boundaries. Therefore, if there are $n_s$ structures each which require a data assimilation node, there must be $n_s - 1$ transition layers which also require a node, and there are always $1$ endpoints. So, using layer-based placement, the minimum number of measurement points should be : 
\begin{align*}
m_h = n_s +n_s- 1 + 1 = 2n_s
\end{align*} data assimilation nodes are required. Using \eqref{mh_est}, we obtain,
\begin{align*}
\lambda & = \frac{L}{n_s} \approx\frac{2L}{m_h} \approx\frac{2L}{c_0\nu^{-p}} \approx \frac{2L}{c_0}\nu^p.
\end{align*}

In this sense, we have found an estimate for $\lambda$ in terms of $\nu$ and $L$. Thus, we have used AOT-style data assimilation to estimate the parameter $\nu$ in terms of measurements of $n_s$. We note that this may be a useful approach to solving certain inverse problems.

\FloatBarrier

 \section*{Acknowledgment}
 \noindent
 The research of A.L. was supported in part by NSF grant no. DMS--1716801.
 The research of C.V. was supported in part by the UCARE program at the University of Nebraska--Lincoln.


\end{document}